\documentclass[11pt]{amsart}
\usepackage{amsmath,amstext,amscd}
\usepackage{amssymb}
\usepackage{epsfig}
\usepackage{amsfonts,latexsym}

\newtheorem{thm}{Theorem}[section]
\newtheorem{lemma}[thm]{Lemma}
\newtheorem{cor}[thm]{Corollary}
\newtheorem{prop}[thm]{Proposition}

\theoremstyle{remark}

\newtheorem{remark}[thm]{Remark}

\newtheorem{Open questions}[thm]{Open questions}

\newtheorem{Problem}[thm]{Problem}
\newtheorem{Open Problem}[thm]{Open Problem}

\newcommand{\cal}{\mathcal}

\newcommand{\F}{\mathbb{F}}
\newcommand{\Z}{\mathbb{Z}}
\newcommand{\N}{\mathbb{N}}

\newcommand{\Diam}{\hbox{\rm Diam}}

\newcommand{\Cay}{\hbox{\rm Cay}}

\newcommand{\R}{\hbox{\rm R}}

\newcommand{\SL}{\hbox{\rm SL}}

\newcommand{\PSL}{\hbox{\rm PSL}}

\newcommand{\ms}{\medskip}

\newcommand{\onto}{{\kern3pt\to\kern-8pt\to\kern3pt}}

\newcommand{\m}{\mathbf{m}}

\renewcommand{\AA}{\cal{A}}
\newcommand{\BB}{\cal{B}}

\newcommand{\KK}{\cal{K}}

\newcommand{\set}[1]{\left\{#1\right\}}

\newcommand{\abs}[1]{\left|#1\right|}
\renewcommand{\ni}{\noindent}
\renewcommand{\ss}{\smallskip}
\renewcommand{\ms}{\medskip}
\newcommand{\bs}{\bigskip}

\newcommand{\Proof}{\ni \emph{Proof. }}

\newcommand{\negsp}{\!\!\!\!\!\!\!\!\!\!\!\!\!\!\!\!\!\!\!\!\!\!\!\!\!\!\!\!\!\!\!\!\!\!\!\!\!\!\!\!\!\!\!}
\newcommand{\possp}{\quad\quad\quad\quad\quad}

\begin{document}

\title{Diameters of Cayley graphs of $\SL_n(\Z/k\Z)$}

\author{M. Kassabov and T. R. Riley}
\thanks{The second author gratefully acknowledges support from NSF grant  0404767.}

\address{Department of Mathematics \\
Cornell University\\
310 Malott Hall\\
Ithaca, NY 14853, USA}
\email{kassabov@aya.yale.edu}
\urladdr{http:/\!/www.math.cornell.edu/\~{}kassabov/}

\address{Mathematics Department,
10 Hillhouse Avenue,
P.O. Box 208283,
New Haven, CT 06520-8283, USA}
\email{tim.riley@yale.edu}
\urladdr{http:/\!/www.math.yale.edu/users/riley/}

\date{January 24, 2005}

\begin{abstract}
We show that for integers $k \geq 2$ and $n\geq 3$, the diameter of the
Cayley graph of $\SL_n(\Z/k\Z)$ associated to a standard two-element
generating set, is at most a constant times $n^2\ln k$.
This answers a question of A.~Lubotzky concerning $\SL_n(\F_p)$ and is unexpected
because these Cayley graphs do not form an expander family.
Our proof amounts to a quick algorithm for finding short words
representing elements of $\SL_n(\Z/k\Z)$.
\end{abstract}

\maketitle

\section{Introduction}

\label{intro}

This paper concerns expressing elements of
$\SL_n(\Z/k\Z)$, for integers $k \geq 2$ and $n \geq 3$, as words in the two-element generating set
$\{\AA_n,\BB_n\}$, where
$$
\AA_n := \left(\begin{array}{ccccc}
1 & 1 &  &  &  \\
  & 1 &  &  &  \\
  &  & 1 &  &  \\
  &  &  & \ddots &  \\
  &  &  &  & 1
\end{array}\right),
\quad
\BB_n :=  \left(\begin{array}{rlllll}
 \ \ \
 0 & 1 &   &        &  \\
   & 0 & 1 &        &  \\
   &   & 0 & \ddots &  \\
   &   &   & \ddots & 1\\
(-1)^{n-1}\!\!\!\!\!\!\!\!\!\!\!\!\!
   &  &    &        & 0
\end{array}\right).
$$

From the point of view of word length, one might suspect this to be an
inefficient generating set because the conjugates of $\AA_n$ by small powers of
$\BB_n$ generate a nilpotent group,  and the diameters of nilpotent groups are large \cite{AB}. However
we show in this paper:

\begin{thm} \label{main thm}
For all integers $k \geq 2$ and $n \geq 3$,
$$
\Diam \ \Cay(\SL_n(\Z/k\Z), \{\AA_n,\BB_n \}) \leq  3600\,n^2 \ln k.
$$
Moreover, there is an algorithm which expresses matrices in $\SL_n(\Z/k\Z)$
as words on $\AA$ and $\BB$ of length $O(\ln \abs{\SL_n(\Z/k\Z)})$ in time $O(\ln \abs{\SL_n(\Z/k\Z)})$.
\end{thm}

The $n^2 \ln k$ term is the best possible because
a logarithm of $\abs{ \SL_n(\Z/k\Z) } \sim k^{n^2 -1}$ gives a
lower bound on the diameter of
$\Cay(\SL_n(\Z/k\Z), \set{\mathcal{A}_n, \mathcal{B}_n})$.
More precise tracking of word length in our arguments would
lead to an improvement of the constant from $3600$ to at least $1400$,
but at the expense of complicating the exposition.

Our result is better than that obtainable by known methods that
use the heavy machinery of Property $T$, Kazhdan constants and
expander families.  For fixed $n \geq 3$, Property \emph{T} of $\SL_n(\Z)$ implies that
$$
\set{\Cay ( \SL_n(\Z/k\Z), \{\AA_n,\BB_n\}) \mid k \geq 2}
$$
is an expander family.  So the diameter of
$\Cay( \SL_n(\Z/k\Z), \{\AA_n,\BB_n\})$
is at most $C(n)\ln k$, where the constant $C(n)$ is
related to the Kazhdan constant $\KK(\SL_n(\Z),\{\AA_n,\BB_n\})$
by
$
C(n) < n^2/\KK(\SL_n(\Z),\{\AA_n,\BB_n\})^2.
$
Lower bounds for Kazhdan constants are hard to come by.
Using the bounds of \cite{Kassabov} for
$\KK(\SL_n(\Z),S)$, where $S$ is the set of all elementary matrices $e_{i,j}$, one can show
$\KK(\SL_n(\Z),\{\AA_n,\BB_n\}) > n^{-3/2}$.
This implies that $C(n) = O(n^5)$.

Were
$$
\set{\Cay ( \SL_n(\Z/k\Z), \{\AA_n,\BB_n\}) \mid k \geq 2, n \geq 3}
$$
an expander family, our $O(n^2 \ln k)$ bound would immediately follow.
But this is not so: on page 105 of~\cite{LW} an argument of
Yael Luz is given that shows the expander constant of
$\Cay ( \SL_n(\Z/k\Z), \{\AA_n,\BB_n\})$
to be at most $5/n$, which is not bounded away from $0$.
(In fact, in~\cite{LW} the generating set considered
has one additional element; none-the-less the argument there
applies to $\{\AA_n,\BB_n\}$.  We remark that it follows that
$\KK(\SL_n(\Z),\{\AA_n,\BB_n\}) \leq \sqrt{2/(n-2)}$.)

Analogous results for $\SL_2(\Z/k\Z)$ and $\SL_2(\F_p)$ cannot be
proved using our methods.  Indeed, there is no known fast
algorithm which writes elements in $\SL_2(\F_p)$ as short words in
$\AA$ and $\BB$. For results in this direction
see~\cite{Dinai, GS, Larsen}.

This article builds on methods in~\cite{Riley5}, where
it is shown (Theorem~5.1) that for all $n \geq 3$, the diameter of $\Cay(\SL_n(\F_p), S)$
is at most a constant times $n^2\ln p$, where $S$ is the set of
all elementary matrices $e_{i,j}$.  By expressing the elementary matrices as
words in $\AA_n$ and $\BB_n$ one deduces~\cite[Corollary~1.1]{Riley5} that the diameter of
$\Cay(\SL_n(\F_p), \{\AA_n,\BB_n\})$ is at most a constant times $n^3 \ln p$.

Theorem~\ref{main thm} affirmatively answers a question of
A.~Lubotzky~\cite[Problem 8.1.3]{Lubotzky} and improves on,
and provides a constructive proof for, a result of~Lubotzky, Babai
and~Kantor:
\begin{prop}[\cite{BKL, Lubotzky}]
\label{BKL prop}
There exists $K >0$ such that for all $n \geq 3$ and primes $p$,
there is a set $\Sigma$ of three generators for $\SL_n(\F_p)$ such that
$$
\Diam \ \Cay(\SL_n(\F_p), \Sigma ) \ \leq  \ K n^2 \ln  p.
$$
\end{prop}

In~\cite{BKL}, it is shown that there is a constant
$K>0$ such that every finite simple non-abelian group $\Gamma$ has
a seven-element generating set $S$ such that
$\Diam \ \Cay(\Gamma,S) \leq K \ln \abs{\Gamma}$.
For $\Gamma=\PSL_{n}(\F_q)$ and $n \geq 10$, Kantor~\cite{Kantor} improved
this by showing that $S$ could be found with only two elements, one of which is an involution.

The methods in this paper can be generalized to show that, with
respect to a generating set consisting of a Weyl element and a
generator of a root subgroup, there exists $K>0$ such that the diameter of
any Chavalley group $\Gamma$ over $\Z/k\Z$, of rank at least $2$,
is at most $K \ln|\Gamma|$.

A further extension shows that for every finite simple group $\Gamma$ of Lie type and of rank at least $2$,
there is a $3$ element generating set $S$ such that
$
\Diam \ \Cay(\Gamma, S )  \leq   K \ln \abs{\Gamma}.
$
Combined with a similar result for the rank $1$ groups from~\cite{BKL} and
the corresponding result for the alternating/symmetric
groups~\cite[Proposition 8.1.6]{Lubotzky}, this yields
\begin{thm}
There exists $K >0$ such that every finite
simple non-abelian group $\Gamma$ has
a $4$ element generating set $S$ such that
$$
\Diam \ \Cay(\Gamma,S) \leq K \ln \abs{\Gamma}.
$$
\end{thm}

Moreover, as all the proofs are suitably constructive, there is a fast algorithm which, given
$g \in \Gamma$ produces a word on $S$ representing $g$, provided that $\Gamma$ is not a factor of a lattice in a rank $1$ Lie group.

\ms

The following problems offer a broader context for the study of
diameters of Cayley graphs of $\SL_n(\Z/k\Z)$ and $\SL_n(\F_p)$.

\begin{Problem}
\label{uniform T}
Fix $n \geq 3$. Does $\SL_n(\Z)$ enjoy \emph{uniform} Property $T$\,?
\end{Problem}

\begin{Problem} \label{uniform diameters}
Fix $n \geq 2$.  Does there exist $K(n)>0$ such that for all
generating sets $X$ for $\SL_n(\Z)$ and all primes $p$
$$\Diam \ \Cay(\SL_n(\Z/k\Z), X)  \ \leq \ K(n) \ln p\,?$$
\end{Problem}

\ni For fixed $n \geq 3$, an affirmative answer to
Problem~\ref{uniform T} would imply an affirmative answer to
 \ref{uniform diameters}.  More details can
be found in~\cite{Lubotzky} and~\cite{LZ}.

\begin{Problem}
What are the maximal and minimal Kazhdan constants for $\SL_n(\Z/k\Z)$ over generating sets of bounded size?
\end{Problem}

\ni
It is shown in~\cite{Kassabov2} that
$ \max \set{ \KK(\SL_n(\Z/k\Z),X) \mid \ \abs{X} <30 } > 10^{-4}$,
and so there are $30$-element generating sets, with respect to which the
Cayley graphs of $\SL_n(\Z/k\Z)$ (with both $n$ and $k$ allowed to vary) form an expander family.

\bs

\section{Generating bit-row and bit-column matrices} \label{bit columns}

All the computations in this and the next section are in $\SL_n(\Z)$ where $n \geq 3$.
Let us fix some notation and terminology.
Denote the matrix with $1$'s on the diagonal and in the $i,j$-th place, and $0$'s
everywhere else by $e_{i,j}$. 
Suppressing their subscripts, we denote $\AA_n$ and $\BB_n$ by $\AA$ and $\BB$.
Define a \emph{row} (\emph{column}) \emph{matrix} to be a square matrix whose diagonal entries are all $1$'s
and which differs from the identity only in one row (column).
A \emph{bit}-row (\emph{bit}-column) matrix is a row (column) matrix whose entries are all in $\set{0,\pm 1}$.
For a sequence $\m = \set{m_i}^n_{i=2}$ define $R_{\m}$ to be the row matrix
whose entries all agree with those of the identity matrix except for those in row $1$ which is
$$
( 1, \  m_2, \ m_3, \ m_4 ,   \ \ldots,  \ m_n).
$$

This section is devoted to proving the following proposition and an analogue concerning column matrices.

\begin{prop} \label{main prop}
Suppose $M \in  \SL_n(\Z)$ is a bit-row matrix. There is a word on
$\AA$ and $\BB$ that represents $M$ and, if the first row of $M$ differs from the identity,  has length at most $48n$, and at most $49n$ otherwise.
\end{prop}

\ms
\Proof
For integers  $s_2, \ldots, s_{n-1}$,
$t_2, \ldots, t_{n-2}$ define
$$
N := \left(\begin{array}{ccccccc}
1 & 0 & 0 &  &  &  &  \\
& 1 & -s_2 & -t_2 &  &  &  \\
&  & 1 & -s_3 & -t_3 &  &  \\
&  &  & 1 & -s_4 & \ddots &  \\
&  &  &  & 1 & \ddots & -t_{n-2} \\
&  &  &  &  & \ddots & -s_{n-1} \\
&  &  &  &  &  & 1
\end{array}\right).
$$
\bs

\begin{lemma} \label{matrix mult}
The matrix $Ne_{1,2}N^{-1}$ is equal to $R_{\m}$, where the sequence $\{m_i\}$ is defined recursively by
$m_2  = 1$, $m_3  =  s_2$ and $m_i  =  m_{i-1}s_{i-1} + m_{i-2} t_{i-2}$ for $4\leq i \leq n$.
\end{lemma}
\ms
\Proof
Rows $2$ to $n$ of $R_{\m}N$ and $Ne_{1,2}$ are the same as rows $2$ to $n$ of $N$. The recursion defining $\set{m_i}$ ensures that the first rows of $R_{\m}N$ and $Ne_{1,2}$ also agree.

\bs

\begin{lemma} \label{solve eqns}
Suppose $m_2=1$ and $m_3, \ldots, m_n \in \set{0,\pm 1,\pm 2}$ satisfy one of the two conditions:
\begin{itemize}
\item[(i)]  $m_{2i} =1$ for all $i$,
\item[(ii)] $m_{2i+1} =1$ for all $i$.
\end{itemize}
Then there exist $s_2, \ldots, s_{n-1} \in \set{0,\pm 1,\pm 2}$ and $t_2, \ldots, t_{n-2} \in \set{0,1}$
satisfying the equations in Lemma~\ref{matrix mult}.
\end{lemma}

\Proof There is a solution in case
(i) with $s_{2i}=m_{2i+1}$, $s_{2i+1}=0$, $t_{2i} =1$, $t_{2i+1}=0$ for all $i \geq 1$,
and in case (ii) with $s_2=1$,  $s_{2i+2}=0$, $s_{2i+1}=m_{2i+2}$, $t_{2i}=0$, $t_{2i+1} =1$, for all $i \geq 1$.

\bs
Our next lemma is an immediate consequence of the previous two.

\begin{lemma} \label{two matrices}
If $m_2=0$ and $m_3, \ldots, m_n \in \set{0,\pm 1}$ then
there exist two matrices $N_1,N_2$ of the same form as $N$,
having all entries in $\set{0,\pm 1, \pm 2}$, and satisfying
$R_{\m} = {N_1}e_{12} {N_1}^{-1} {N_2}{e_{12}}^{-1} {N_2}^{-1}$.
\end{lemma}

\ms

\begin{lemma}
\label{word for N}
If $s_2, \ldots, s_{n-1} \in \set{0,\pm 1,\pm 2}$ and $t_2, \ldots, t_{n-2} \in \set{0,\pm 1}$
then the word
$$
w \ := \ \BB^{-(n-2)}Q_{n-1} \BB Q_{n-2} \BB  \ldots Q_2 \BB,
$$
where $Q_i := \AA^{-s_i} [\AA,\BB^{-1}\AA\BB]^{-t_i}$ and $t_{n-1}:=0$,
represents $N$ and has length at most $12n-24$ as a word on $\AA$ and $\BB$.
\end{lemma}

\ms

\Proof This result follows from the observations that $w$ equals
\begin{equation*}
(\BB^{-(n-2)} Q_{n-1} \BB^{n-2})
(\BB^{-(n-3)} Q_{n-2} \BB^{n-3})
 \ldots (\BB^{-1} Q_2 \BB),
\end{equation*}
 in $\SL_n(\Z)$,
and $(\BB^{-(i-1)} Q_i \BB^{i-1})$ equals the row matrix
whose entries in the $i$-th row agree with those of $N$ and
whose remaining entries agree with the identity matrix.

\bs

To complete the proof of Proposition~\ref{main prop} in the row matrix case,
first suppose $M = R_{\m}$ where $m_2=0$ and
$m_3, \ldots, m_n \in \set{0,\pm 1}$.
Lemmas~\ref{two matrices} and~\ref{word for N} supply a word
$w_{\m}$ on $\AA$ and $\BB$ that represents $M$ and has length at most $4(12 n - 24)+2$.
We can change $m_2$ to $\pm 1$ by right-multiplying by ${e_{1,2}}^{\pm 1}= \AA^{\pm 1}$.
This proves that $\ell(w_{\m}) < 48n$, as claimed.

Conjugating a matrix by a power of $\BB$ moves its entries diagonally
(changing some of their signs if $n$ is even).
So any given bit-row matrix $M$ equals $\BB^l R_{\m} \BB^{-l}$ for some $\m=\set{m_i}_{i=2}^n$ with $m_i\in \set{0,\pm1}$, and some $- n/2 \leq l \leq n/2$.  The cost to word length of conjugating is at most $n$ and so
$M$ can be written as a word on $\AA$ and $\BB$ of length at most $49n$. \qed

\bs
There is a natural analog of Proposition~\ref{main prop} for bit-column matrices:
\begin{prop} \label{main prop_col}
If $M \in  \SL_n(\Z)$ is a bit-column  matrix then there is a word on
$\AA$ and $\BB$ that represents $M$ and has length at most $48n$
if the final column differs from the identity, and at most $49n$ otherwise.
\end{prop}

\Proof As in Proposition~\ref{matrix mult}, if we define
$$
\widetilde{N} := \left(\begin{array}{ccccccc}
1 & -s_1 & -t_1 &  &  &  &  \\
& 1 & -s_2 & -t_2 &  &  &  \\
&  & \ddots & \ddots & \ddots &  &  \\
&  &  & 1 & -s_{n-3} & -t_{n-3} &  \\
&  &  &  & 1 & -s_{n-2} & 0 \\
&  &  &  &  & 1 & 0 \\
&  &  &  &  &  & 1
\end{array}\right),
$$
we find $\widetilde{N}^{-1}e_{n-1,n}\widetilde{N}$ is a column matrix $C_{\m'}$ whose final column is
$$
(m'_1, \, \ldots, m'_{n-1}, \, 1)^{\textit{tr}},$$
where
$m'_{n-1}  = 1$, $m'_{n-2}  =  s_{n-2}$ and $m'_i  =  m'_{i+1}s_{i} + m'_{i+2} t_{i}$
for $1\leq i \leq n-3$.

The analogue of Lemma~\ref{solve eqns},
implies that if $m'_1, \ldots, m'_{n-2} \in \set{0,\pm 1,\pm 2}$
satisfies $m'_{2i}=1$ for all $i$, or $m'_{2i+1}=1$ for all $i$, then
$C_{\m'}$ can be written as a short word on $\AA$ and $\BB$.
And, as in Lemma~\ref{two matrices}, two such column matrices can be multiplied to
produce any given bit-column matrix $C_{\m''}$ in which
$m''_{n-1}=0$.  Then $m''_{n-1}$ can then be changed to $\pm 1$ if required by left-multiplying by  ${e_{n-1,n}}^{\pm 1}$.

In this context, the analogue of the length bound of Lemma~\ref{word for N} is $12n - 26$.
Therefore the cost of producing $C_{\m''}$ is at most  $4(12 n - 26)+ 10$ because
$e_{n-1,n} = \BB^{2} \AA \BB^{-2}$ if $n$ is odd and
$e_{n-1,n} = \BB^{2} \AA^{-1} \BB^{-2}$ if  $n$ is even, which has word length $5$ in both cases.
This also shows that the cost of altering  $m''_{n-1}$ is at most
$5$. The total cost, then, is within the claimed bound of $48n$.
As in the row matrix case, it follows that every bit-column matrix can
be written as a word on $\AA$ and $\BB$ of length at most $49n$.  \qed

\begin{remark}
This construction yields an algorithm with running time $O(n)$ which
produces a short word on $\AA$ and $\BB$ representing any given
bit-row or bit-column matrix in $\SL_n(\Z)$.
\end{remark}
\begin{remark}
If we allow $s_i\in \set{1,2,4}$ in $N$ and we set all the $t_i=0$
(so $N$ has one super-diagonal, not two) then it is possible to find shorter words representing
certain row and column matrices in $\SL_n(\Z)$ whose entries are particular powers of $2$.   If $k$ is odd we can use these row and column matrices to
obtain a better upper bound than that of Theorem~\ref{main thm}.
This breaks down when $k$ is even because there are insufficient
invertible elements in $\F_2$.
\end{remark}

\section{Generating row and column matrices}

Before we come to the main result of this section we give a lemma which is
essentially~\cite[Lemma~2.2]{Riley5}. It concerns expressing matrices
${e_{i,j}}^{F_{2l}}$ and ${e_{i,j}}^{F_{2l+1}}$, where the powers are Fibonacci numbers
(defined recursively by $F_0=0$,  $F_1=1$, and  $F_{i+2}=F_{i+1}+F_i$),
as short words on $\set{e_{i,j} \mid i \neq j}$.
This lemma will be superseded by Lemma~\ref{compress}, but the detailed
calculation we give in the proof of this simpler case is key to understanding
the proof of the more general result.

\begin{lemma}
\label{Fib powers}
For  non-negative integers $l$, the words
$$
\begin{array}{cc}
{e_{1,2}}^2
(e_{2,1}e_{1,2})^l e_{1,3}(e_{2,1}e_{1,2})^{-l}
{e_{1,2}}^{-1}
(e_{2,1}e_{1,2})^{l}{e_{1,3}}^{-1}(e_{2,1}e_{1,2})^{-l}
{e_{1,2}}^{-1}, &
\mbox{ and } \\
{e_{1,2}}^2
(e_{2,1}e_{1,2})^l e_{2,3} (e_{2,1}e_{1,2})^{-l}
{e_{1,2}}^{-1}
(e_{2,1}e_{1,2})^{l} {e_{2,3}}^{-1} (e_{2,1}e_{1,2})^{-l}
{e_{1,2}}^{-1}
\end{array}
$$
equal ${e_{1,3}}^{F_{2l}}$ and ${e_{1,3}}^{F_{2l+1}}$, respectively,
in $\SL_3(\Z)$.
\end{lemma}
\Proof
We multiply out the first of these words from left to right
as follows.  The calculation for the second is similar.  The
notation for each step shown is $S\stackrel{T}{\longrightarrow} ST$.
\ms
$$
\begin{array}{l@{}l@{}l}
\left(\begin{array}{ccc}
1 & 0 & 0 \\
0 & 1 & 0 \\
0 & 0 & 1
\end{array}\right)
\stackrel{{e_{1,2}}^2}{\xrightarrow{\quad\quad}}
\left(\begin{array}{ccc}
1 & 2 & 0 \\
0 & 1 & 0 \\
0 & 0 & 1
\end{array}\right)
&
\stackrel{(e_{2,1}e_{1,2})^l}{\xrightarrow{\quad\quad}}
&
\left(\begin{array}{ccc}
F_{2l+2} & F_{2l+3} & 0 \\
F_{2l} & F_{2l+1} & 0  \\
0 &  0 &  1
\end{array}\right)
\\
\quad\quad
\stackrel{{e_{1,3}}}{\xrightarrow{\quad\quad}}
\rule{0pt}{12mm}\left(\begin{array}{ccc}
F_{2l+2} & F_{2l+3} & F_{2l+2}\\
F_{2l} & F_{2l+1} & F_{2l} \\
0 & 0 & 1
\end{array}\right)
&
\stackrel{(e_{2,1}e_{1,2})^{-l}}{\xrightarrow{\quad\quad}}
&
\left(\begin{array}{ccc}
1 & 2 & F_{2l+2} \\
0 & 1 & F_{2l} \\
0 & 0 & 1
\end{array}\right)
\\
\quad\quad
\stackrel{{e_{1,2}}^{-1}}{\xrightarrow{\quad\quad}}
\rule{0pt}{12mm}\left(\begin{array}{ccc}
1 & 1  & F_{2l+2}\\
0 & 1 & F_{2l} \\
0 & 0 & 1
\end{array}\right)
&
\stackrel{(e_{2,1}e_{1,2})^l}{\xrightarrow{\quad\quad}}
&
\left(\begin{array}{ccc}
F_{2l+1} & F_{2l+2} & F_{2l+2} \\
F_{2l} & F_{2l+1} & F_{2l}  \\
0 & 0  & 1
\end{array}\right)
\\
\quad\quad
\stackrel{{e_{1,3}}^{-1}}{\xrightarrow{\quad\quad}}
\rule{0pt}{12mm}\left(\begin{array}{ccc}
F_{2l+1} & F_{2l+2} & F_{2l} \\
F_{2l} & F_{2l+1} & 0   \\
0 & 0 & 1
\end{array}\right)
&
\stackrel{(e_{2,1}e_{1,2})^{-l}}{\xrightarrow{\quad\quad}}
&
\left(\begin{array}{ccc}
1 & 1 & F_{2l} \\
0 & 1 & 0 \\
0 & 0 & 1
\end{array}\right)
\\
\quad\quad
\stackrel{{e_{1,2}}^{-1}}{\xrightarrow{\quad\quad}}
\rule{0pt}{12mm}\left(\begin{array}{ccc}
1 & 0 & F_{2l} \\
0 & 1 & 0  \\
0 & 0 & 1
\end{array}\right)
& & \hfill\qed
\end{array}
$$

\bs

\begin{prop} \label{main prop2}
Suppose $M \in  \SL_n(\Z)$ is a row or column matrix with entries in
$\set{-K+1,\ldots, 0, \ldots, K-1}$, where $K \geq 1$.
Then there is a word on $\AA$ and $\BB$, representing $M$, that has length at
most $1200\,n \ln K + 400n$.
\end{prop}

\Proof
The proof in the row matrix case generalizes Lemma~\ref{Fib powers} -- instead of
using $e_{1,3}$ and $e_{2,3}$ we use general bit-row matrices; they allow the simultaneous
construction of sums of Fibonacci numbers in entries $3, \ldots, n$ of the first row.
These sums of Fibonacci numbers are as per Zeckendorf's
Theorem~\cite{GKP, Zeckendorf}, which  states that every nonzero integer $m$ can be expressed in a unique way as
\begin{equation*} \label{Zeckendorf}
m= \pm(F_{l_1} + F_{l_2}+ \cdots + F_{l_r}),
\end{equation*}
with $l_1\geq 2$ and $l_{j+1} - l_j \geq 2$ for all $1 \leq j < r$.
This result can be proved by an easy induction argument and, in fact,
$F_{l_r}$ is the largest Fibonacci number no bigger
than $\abs{m}$, and $F_{l_{r-1}}$ is the largest no bigger than
$\abs{m}- F_{l_r}$, and so on.  Since
$F_s = (\tau^s - (-\tau)^{-s})/\sqrt{5}$ for all $s \in \N$,
where $\tau := (1+\sqrt 5)/2$,
we get
$F_s \geq (\tau^s -1)/\sqrt{5}$. Thus, as $F_{l_r} \leq \abs{m}$, we find
$$
l_r  \leq  \log_{\tau} (1+ \abs{m} \sqrt 5 ) < 2+3 \ln \abs{m},
$$
from which we derive the bound on $L$ in the following lemma.

\begin{lemma} \label{compress}
Suppose $\m:= \set{m_i}^n_{i=3}$ is a sequence of integers, such that $|m_i| < K$ for all $i$.
As per Zeckendorf's Theorem, write
$$
m_i \ = \   \sum_{j=1}^{L} (c_{ij} F_{2j} +  d_{ij} F_{2j+1})
$$
where $c_{ij},d_{ij} \in \set{0,\pm 1}$ and $L \leq (2+3 \ln K)/2 - 1/2$.
Let $u_{\m}$ be the word
$$(e_{2,1}e_{1,2})a_{1}b_{1} (e_{2,1} e_{1,2}) a_2b_2(e_{2,1}e_{1,2})\ldots (e_{2,1}e_{1,2})a_{L}b_{L}$$
in which
$a_j$ is the row matrix with first row $(1, 0, c_{3j}, \ldots, c_{nj})$ and $b_j$ is
the row matrix with second row $(0, 1, d_{3j}, \ldots, d_{nj})$.
Let  $v_{\m}$ be the word obtained from $u_{\m}$ by replacing each $a_j$ and $b_j$
by its inverse. Define
$$
w_{\m} :=   {e_{1,2}}^2u_{\m} (e_{2,1}e_{1,2})^{-L}{e_{1,2}}^{-1}v_{\m}(e_{2,1}e_{1,2})^{-L} {e_{1,2}}^{-1}.
$$
Then in $\SL_n(\Z)$ the row matrix with first row $(1,0,m_3, \ldots, m_n)$ is represented by $w_{\m}$.
\end{lemma}
\Proof
Lemma~\ref{Fib powers} gives the special cases of this lemma in which
$n=3$ and $m_3$ is $F_{2l}$ or $F_{2l+1}$. Below we multiply out $w_{\m}$ from left to right,
using a more general and concise version of the calculation used to prove
Lemma~\ref{Fib powers}.  We display the top two rows only; all others agree with the
identity matrix throughout the calculation.  All the summations range over $j = 1, \ldots, L$.
$$
\begin{array}{c@{\!\!\!}c@{\!}c@{\!}c@{\!\!}c@{\!\!}c@{\!\!\!}c}
\Bigg(
&
\begin{array}{c}1\\0\end{array}
&
\begin{array}{c}2\\1\end{array}
&
\begin{array}{c}0\\0\end{array}
&
\begin{array}{c}\cdots\\ \cdots\end{array}
&
\begin{array}{c}0\\0\end{array}
&
\Bigg)
\\
& & & \possp \rule{0mm}{8mm}\Big\downarrow \parbox{0mm}{$\scriptstyle u_{\m}$\negsp} & & &
\\
\rule{0mm}{10mm}\Bigg(
&
\begin{array}{c}F_{2L+2}\\F_{2L}\end{array}
&
\begin{array}{c}F_{2L+3}\\F_{2L+1}\end{array}
&
\begin{array}{c}\sum (c_{3j} F_{2j+2} +  d_{3j} F_{2j+3}) \\ \sum(c_{3j} F_{2j} +  d_{3j} F_{2j+1})\end{array}
&
\begin{array}{c}\cdots\\ \cdots\end{array}
&
\begin{array}{c}\sum (c_{nj} F_{2j+2} +  d_{nj} F_{2j+3}) \\ \sum(c_{nj} F_{2j} +  d_{nj} F_{2j+1})\end{array}
&
\Bigg)
\\
& & & \possp \rule{0mm}{8mm}\Big\downarrow \parbox{0mm}{ $\scriptstyle (e_{2,1}e_{1,2})^{-L}$\negsp} & & &
\\
\rule{0mm}{10mm}\Bigg(
&
\begin{array}{c}1\\0\end{array}
&
\begin{array}{c}2\\1\end{array}
&
\begin{array}{c}\sum(c_{3j} F_{2j+2} +  d_{3j} F_{2j+3}) \\ \sum(c_{3j} F_{2j} +  d_{3j} F_{2j+1})\end{array}
&
\begin{array}{c}\cdots\\ \cdots\end{array}
&
\begin{array}{c}\sum(c_{nj} F_{2j+2} +  d_{nj} F_{2j+3}) \\ \sum(c_{nj} F_{2j} +  d_{nj} F_{2j+1})\end{array}
&
\Bigg)
\\
& & & \possp \rule{0mm}{8mm}\Big\downarrow \parbox{0mm}{$\scriptstyle {e_{1,2}}^{-1}v_{\m}$\negsp}& & &
\\
\rule{0mm}{10mm}\Bigg(
&
\begin{array}{c}F_{2L+1}\\F_{2L}\end{array}
&
\begin{array}{c}F_{2L+2}\\F_{2L+1}\end{array}
&
\begin{array}{c}\sum(c_{3j} F_{2j} +  d_{3j} F_{2j+1}) \\ 0 \end{array}
&
\begin{array}{c}\cdots\\ \cdots\end{array}
&
\begin{array}{c}\sum(c_{nj} F_{2j} +  d_{nj} F_{2j+1}) \\ 0 \end{array}
&
\Bigg)
\\
& & & \possp \rule{0mm}{8mm}\Big\downarrow \parbox{0mm}{$\scriptstyle (e_{2,1}e_{1,2})^{-L}{e_{1,2}}^{-1}$\negsp} & & &
\\
\rule{0mm}{10mm}\Bigg(
&
\begin{array}{c}1 \\ 0 \end{array}
&
\begin{array}{c}0 \\ 1 \end{array}
&
\begin{array}{c}\sum(c_{3j} F_{2j} +  d_{3j} F_{2j+1}) \\ 0 \end{array}
&
\begin{array}{c}\cdots\\ \cdots\end{array}
&
\begin{array}{c}\sum(c_{nj} F_{2j} +  d_{nj} F_{2j+1}) \\ 0 \end{array}
&
\Bigg
).
\end{array}
$$
\ss

\ni The sums in this final matrix are, by definition, equal to
 $m_3, \ldots, m_n$ and so the lemma is proved.

\bs

Returning to the proof of Proposition~\ref{main prop2}, note that a
conjugate of $M$ by a power of $\BB$ is a row matrix $\R_{\m}$
in which the first row is $(1, m_2, m_3, \ldots, m_n)$.
On the alphabet $\AA$ and $\BB$, we find $e_{1,2}=\AA$ and
so has length $1$, and $e_{2,1}, a_{j},b_j$
are all bit-row matrices and so, by Proposition~\ref{main prop},
can be expressed as words of length at most $48n$.
So the word $w_{\m}$ of Lemma~\ref{compress} can be re-expressed as a word on $\AA$ and $\BB$
of length $390nL$, where the contributions to this estimate are
$$
\begin{array}{r@{\,\,\,\times\,\,\,}c@{\;\;\mbox{from}\;\;}l}
4+4L &  1  & e_{1,2} \\
4L   & 48n & e_{2,1} \\
2L   & 48n & a_j \\
2L   & 48n & b_j.
\end{array}
$$

A revised version of Lemma~\ref{compress} in which we build up Fibonacci numbers in columns
$1$ and $n$ using $e_{1,n}$ and $e_{n,1}$ rather than in columns $1$ and  $2$ using $e_{1,2}$ and $e_{2,1}$,
produces a word on $\AA$ and $\BB$ that represents
the row matrix with first row $(1,m_2, 0, \ldots, 0)$.  Mildly revising
the estimates above, we check that the length of this word is at most $390nL$.
Multiplying the two words together gives a word of length at most $780nL$
that represents $R_{\m}$.  Conjugating by a power of $\BB$ recovers $M$
at a further expense to word length of at most $n$.  Then, using the bound on $L$ in Lemma~\ref{compress},
we learn that $M$ can be represented by a word on $\AA$ and $\BB$ of length at most $1200  n \ln K + 400n$.

Obtain the same bound in the column matrix case by transposing and using Proposition~\ref{main prop_col} in place of \ref{main prop}:
reverse the orders of the terms in $w_{\m}, u_{\m}$ and $v_{\m}$, interchange the
$e_{1,2}$'s and $e_{2,1}$'s, and make the $a_i$ and $b_i$ bit-column
matrices rather than bit-row matrices.
\qed

\begin{remark}
It follows from the construction above that there is an algorithm
with running time $O(n\ln K)$ which
produces a short word in $\AA$ and $\BB$
that represents any given row matrix in $\SL_n(\Z)$ with entries in $\{-K+1, \dots, 0,\dots, K-1\}$.
\end{remark}

\bs

\section{The diameter of $\SL_n(\Z/k\Z)$.}
\label{diam}

\ni \textit{Proof of Theorem~\ref{main thm}.}
All row matrices in $\SL_n(\Z/k\Z)$ come from
row matrices in $\SL_n(\Z)$ with entries of absolute value less than $k/2$ and so
can be represented by short words on $\AA$ and $\BB$ as per Proposition~\ref{main prop2}.
So Lemma~\ref{3n matrices} below completes the proof of the bound in Theorem~\ref{main thm}.

Our proof is constructive and amounts to an algorithm for expressing matrices in
$\SL_n(\Z/k\Z)$ as words on $\AA_n$ and $\BB_n$ with running time
$$
O(n^2 \ln k)   =   O(\ln\abs{ \SL_n(\Z/k\Z) }),
$$
provided that $k$ is decomposed as
a product of prime numbers.\qed

\ms
We start with a technical lemma which is also valid for rings satisfying the Bass
stable range condition -- see \cite{Va}.

\begin{lemma} \label{lemm}
Let $a,b\in \Z/k\Z$.  Then there exists $s \in \Z/k\Z$ such that the ideal generated by $a$ and $b$ is the same
as the ideal generated by $a+sb$.
\end{lemma}
\Proof
If $k= \prod {p_i}^{m_i}$ then
$$
\Z/k\Z \simeq \prod \Z/{p_i}^{m_i}\Z
$$
 by the Chinese Remainder Theorem.
Let $a_i$ and $b_i$ be the components of $a$ and $b$ in $\Z/{p_i}^{m_i}\Z$. Define
$s_i:=0$ if the ideal generated by $a_i$ in $\Z/{p_i}^{m_i}\Z$ contains $b_i$
and $s_i:=1$ otherwise. Let $s$ be the element in $\Z/k\Z$ with components $s_i$.
By construction, the components of $a+sb$ are $a_i+s_ib_i$ and in the ring $\Z/{p_i}^{m_i}\Z$
the ideal generated by $a_i+s_ib_i$ is the same as the ideal generated by $a_i$ and $b_i$.
\qed

\begin{cor}
\label{co}
Suppose $\{a_i\}_{i=1}^l$ are elements of $\Z/k\Z$ such that the ideal they
generate is the whole ring. Then there exist $\{t_i\}_{i=2}^l$ such that
$$
a_1 + t_2 a_2 + t_3a_3 + \dots + t_l a_l
$$
is invertible in $\Z/k\Z$.
\end{cor}

In fact, (given the decomposition of $k$ into prime factors)
we can write a fast algorithm to find these coefficients.
This is because $a_i$ and $b_i$ of Lemma~\ref{lemm} can be found
quickly, being $a \, \textup{mod} \,  {p_i}^{m_i}$ and $b \, \textup{mod} \, {p_i}^{m_i}$, respectively.
The maximal power of $p_i$ dividing $a_i$ and $b_i$ determines $s_i$. And in the proof of Lemma~\ref{lemm} we can use $k\sum s_i/{p_i}^{m_i}$, which is easier to compute.

\begin{lemma} \label{3n matrices}
If  $M \in \SL_n(\Z/k\Z)$ then the matrix $M$ can be written as
a product of $n$ row matrices, $n$ column matrices and $n$ elementary matrices.
\end{lemma}

\Proof
We use a version of Gauss-Jordan elimination to prove by induction on
$r=0, \ldots, n$ that $M$ can be transformed to a matrix in which the top
$r$ rows agree with the identity matrix by left- and right-multiplying by a
total of $3r$ row, column and elementary matrices.

The base step $r=0$ holds vacuously.  For the induction step assume $r <n$ and
the top $r$ rows agree with the identity matrix.  If the
final entry on the $(r+1)$-st row is not invertible in $\Z/k\Z$ then, using Corollary~\ref{co},
it can be made invertible by right-multiplying by some column matrix,
because the ideal generated by the $(r+1)$-st to $n$-th entries in row $r+1$ is the whole ring $\Z/k\Z$.
Then make the $r+1,r+1$-entry $1$ by right-multiplying by the appropriate power of $e_{n,r+1}$.
 Then clear all the off-diagonal entries in row $r+1$ by
right-multiplying by the appropriate row matrix.   \qed

\begin{remark}
The constructions in this paper can be used to express matrices $M \in \SL_n(\Z)$ as short words on
$\AA$ and $\BB$ (cf.\ \cite[Theorem~4.1]{Riley5}). However, the resulting upper bounds on word length
are not very good because if we express $M$ as a product of row matrices $R_i$ then the absolute
values of the entries in the $R_i$ may be significantly larger than the absolute values of the entries in $M$.
\end{remark}

\bibliographystyle{plain}
\bibliography{bibli}

\end{document}